\documentclass[a4paper,11pt,reqno]{amsart}
\usepackage{latexsym}
\usepackage{csquotes}
\usepackage{amsmath,amsthm,amssymb,amscd}
\usepackage{fullpage}
\usepackage{xcolor}
\usepackage{stmaryrd}
\usepackage{parskip}
\usepackage{mathtools}
\mathtoolsset{showonlyrefs}
\usepackage{tikz}
\usetikzlibrary{positioning}
\usepackage[colorlinks=false,urlbordercolor=white, pdfencoding=unicode,pdftex,bookmarks=tr
ue,bookmarksnumbered]{hyperref}
  \def\gn#1#2{{$\href{http://groupnames.org/\#?#1}{#2}$}}
\def\gn#1#2{$#2$}  % comment this line out to get html links
\tikzset{sgplattice/.style={inner sep=1pt,norm/.style={red!50!blue},char/.style={blue!50!black},
  lin/.style={black!50}},cnj/.style={black!50,yshift=-2.5pt,left=-1pt of #1,scale=0.5,fill=white}}

\usepackage[activate={true,nocompatibility},final,tracking=true,kerning=true,spacing=true]{microtype}
\usepackage[pdfencoding=unicode,pdftex,bookmarks=true,bookmarksnumbered]{hyperref}
\definecolor{citecolour}{rgb}{0.0, 0.0, 0.8}
\definecolor{urlcolour}{rgb}{1,0.5,0}
\colorlet{linkcolour}{green!50!black}
\hypersetup{colorlinks,breaklinks,
		linkcolor=linkcolour,citecolor=citecolour,
		filecolor=linkcolour, urlcolor=urlcolour}
\usepackage{version}

\theoremstyle{plain}

\newtheorem{theorem}{Theorem}[section]
\newtheorem{proposition}[theorem]{Proposition}
\newtheorem{lemma}[theorem]{Lemma}
\newtheorem{corollary}[theorem]{Corollary}

\theoremstyle{definition}
\newtheorem*{ack}{Acknowledgements}

\newtheorem{definition}[theorem]{Definition}
\newtheorem{remark}[theorem]{Remark}

\numberwithin{equation}{section}
\usepackage[shortlabels]{enumitem}
\begingroup
 \makeatletter
 \@for\theoremstyle:=definition,remark,plain\do{%
  \expandafter\g@addto@macro\csname th@\theoremstyle\endcsname{%
  \addtolength\thm@preskip\parskip
  }%
  }

\DeclareMathOperator{\X}{\mathcal{X}}
\DeclareMathOperator{\W}{\mathcal{W}}
\DeclareMathOperator{\B}{\mathcal{B}}

\newcommand{\R}{\mathbb{\R}}

\newcommand{\pc}[1]{{\mathcal{P}}_{#1}}
\renewcommand{\leq}{\leqslant}
\renewcommand{\geq}{\geqslant}

\newenvironment{proofof}{{\bf {Proof.} }}{\hfill $\blacksquare$ \\}
\newenvironment{proofofmain}{{\bf {Proof of Theorem~\ref{Thm:main}.}}}{\hfill $\blacksquare$ \\}
\begin{document}

\title{The $\X$-series of a  $p$-group and complements of abelian subgroups}

\author{Stefanos Aivazidis}
\address{Department of Mathematics \& Applied Mathematics, University of Crete, Greece}
\email{s.aivazidis@uoc.gr}

\author{Maria Loukaki}
\address{Department of Mathematics \& Applied Mathematics, University of Crete, Greece}
\email{mloukaki@uoc.gr}

\thanks{The first author is partially supported by the Hellenic Foundation for Research and Innovation, Project HFRI-FM17-1733.}

\begin{abstract}
Let $G$ be a  $p$-group. 
We denote by $\X_i(G)$ the intersection of all subgroups of $G$ 
having index $p^i$, for $i \leq \log_p(|G|)$. 
In this paper, the newly introduced series $\{\X_i(G)\}_i$ 
is investigated and a number of results 
concerning its behaviour are proved. 
As an application of these results, 
we show that if an abelian subgroup $A$ of $G$ 
intersects each one of the subgroups $\X_i(G)$ at $\X_i(A)$, 
then $A$ has a complement in $G$. Conversely if an arbitrary subgroup $H$ of $G$ has a normal complement, then $\X_i(H) = \X_i(G) \cap H$.
\end{abstract}



\maketitle

%\tableofcontents

\section{Introduction}\label{Sec:1}
In this paper, we shall only be concerned with
finite groups.
We begin with some general remarks on complementation.
For a fuller discussion in textbook format,
the reader is invited to consult the recent
monograph by Kirtland~\cite{kirtland}.

The (rather broad) question we ask is the following:
\begin{quote}
Let $G$ be a  group and let $N$ be a
normal subgroup of $G$. Under which circumstances is $N$ guaranteed to have
a complement in $G$?
\end{quote}
Perhaps the most well-known and most general
answer is the celebrated Schur-Zassenhaus
theorem which asserts that if $N$ is a Hall
subgroup of $G$ (so that $\gcd(|N|,|G:N|) = 1$), 
then $N$ does indeed have a complement.
In fact, the
theorem asserts more: all complements of $N$ in $G$ are $G$-conjugate. 
The conjugacy part of the theorem 
relies on a solubility assumption: 
either $N$ or $G/N$ must be soluble for it to admit a proof. 
This, in turn, can only be made to hold unconditionally by an appeal 
to the deep Odd Order Theorem of Feit and Thompson~\cite{FTOdd}.

Another answer to the question above
(relaxing the requirement in the abelian case
of the Schur-Zassenhaus theorem)
is furnished by a theorem of Gasch\"{u}tz~\cite{Gaschutz:1952}.
For recent work discussing Gasch\"{u}tz's theorem,
see~\cite{sambale}.

\begin{theorem}[Gasch\"{u}tz]
Let $N$ be an abelian normal subgroup of a  group $G$. Let $N \leq H \leq G$ such that $N$ has a complement in $H$ and $\gcd(|N|,|G : H|) = 1$. 
Then $N$ has a complement in $G$.
\end{theorem}

We see therefore that if $N$ is normal in $G$ and abelian
then $N$ has a complement in $G$ if and only if for each
prime $p$ dividing $|N|$ the unique Sylow $p$-subgroup of $N$ (which is normal in $G$)
has a complement in some (and thus in every) Sylow $p$-subgroup of $G$. 
In fact, the previous statement
is a Reduktionsatz in Gasch\"{u}tz's paper.

This global-to-local reduction that Gasch\"{u}tz's theorem achieves 
reaches an obvious and natural limit 
if $G$ is itself a $p$-group. 
Thus we are led to ask:

\begin{quote}
Let $G$ be a  $p$-group 
and $N$ be a (possibly normal) 
abelian subgroup of $G$. 
Is there a necessary and sufficient condition 
for $N$ to have a complement in $G$?
\end{quote}

In fact, the previous question was 
the exact point of departure for the work reported here.

Although we have not been able to formulate such a condition,
we have found a condition which is sufficient (and in some cases necessary).
We delay presenting the said result
to discuss properties of the $\X$-series
of a  $p$-group first.
This is done in Section~\ref{Sec:2}.
Briefly, given a  $p$-group $G$,
we define $\X_i(G)$ to be 
the intersection of all subgroups of $G$ 
having index $p^i$ in $G$ 
for $i \leq \log_p(|G|)$
and (for technical reasons) extend its definition
to all non-negative integers.
We call the associated series $\{\X_i(G)\}_i$ 
the $\X$-series of $G$. 
This series has several nice properties:
\begin{itemize}
\item It is subgroup-monotone; cf. Proposition~\ref{Prop:SubgMonotone}.
\item It respects direct products; cf. Theorem~\ref{Thm:DirectProd}.
\item For every $i, j \geq 0$ we have  $\X_i (\X_j(G) ) \leq \X_{i+j} (G)$.
Moreover, successive quotients 
of the $\X$-series of $G$ are elementary
abelian $p$-groups; cf. Lemma~\ref{Lem:LikeLCS}.
\item For every $i \geq 0$ we have $\gamma_{p^{i-1} + 1}(G) \leq \X_i(G)$; cf. Proposition ~\ref{Prop:GammaInX}.   
\end{itemize}

In Section~\ref{Sec:3} we resume our discussion
by focussing first on  abelian groups and
what is known in that setting as regards complementation.
We then utilise the tools assembled in Section~\ref{Sec:2} 
to prove our main result:
\begin{theorem}\label{Thm:main}
Let $A$ be an abelian subgroup of the 
$p$-group $G$ and suppose that 
\[
A \cap \X_i(G) = \X_i(A) 
\]
for all indices $i \geq 0$. 
Then $A$ has a complement in $G$. 
\end{theorem}

The importance of this result and its relevance to other
known theorems will be further elucidated in due course.
A sort of partial converse of
Theorem~\ref{Thm:main} holds and its proof is given in Section~\ref{Sec:2}:
\begin{theorem}\label{Thm:converse}
Let $H$ be an arbitrary  subgroup of the $p$-group $G$ and assume that $H$  has a normal complement in $G$. Then 
\[
H \cap \X_i(G) = \X_i(H) 
\]
for all indices $i \geq 0$. 
\end{theorem}

\subsection{Notation}
We outline below some notational
conventions that we will use throughout the paper.
\begin{enumerate}[label={\upshape(\roman*)}]
\item $[n]$ denotes the set of the first $n$ positive integers
$\{1,2,\ldots,n\}$.

\item Let $G$ be a  $p$-group and let $i$
be a non-negative
integer. We will write 
\[
\mho_i(G) = \big\langle g^{p^i} : g \in G \big\rangle
\]
for the subgroup generated by the $p^i$-powers
of elements of $G$. 
Clearly, $\mho_0(G) = G$. The index $i=1$ is sometimes
dropped, so that $\mho_1(G) = \mho(G)$.
These subgroups, for various 
indices $i$, are referred to as the agemo subgroups of $G$. 
We remark that $G^{p^i}$ is sometimes used instead of $\mho_i(G)$ in the literature.

\item We will have occasion to refer to the Frattini
series of a group. For an arbitrary  group,
the Frattini series is defined by $\Phi_0(G) = G$
and $\Phi_i(G) = \Phi(\Phi_{i-1}(G))$,
where $\Phi(G)$ is the Frattini subgroup of $G$,
i.e. the intersection of all maximal subgroups of $G$.
\end{enumerate}

\section{Properties of the $\X$-series}\label{Sec:2}
We begin in earnest by defining the $\X$-series
of a  $p$-group.

\begin{definition}
Given a  $p$-group $G$ and a positive integer $i \leq \log_p |G|$ we write 
\[
\pc{i}= \pc{i}(G) \coloneqq 
\left\{T \leq G \, : \, |G:T| = p^i\right\}
\]
for the collection of subgroups of $G$ of index $p^i$ in $G$. 
We also write 
\[
\X_i(G) \coloneqq \bigcap_{T \in \pc{i}(G)} T\,;
\]
that is, $\X_i(G)$ is the intersection 
of all subgroups of $G$ that  have index $p^i$ in $G$ for all such $i$. 
We call the associated series the $\X$-series of $G$ 
and by convention we write $\X_0(G) = G$.
Finally, if 
$j > \log_p|G|$ then $\X_j(G) = 1$ by assumption.
\end{definition}

Some first properties following the definition of the $\X$-series 
are collected in the following lemma.

\begin{lemma}\label{Lem:BasicProperties}
If $G$ is a  $p$-group then 
the subgroups $\X_i(G)$ are characteristic subgroups of $G$ and 
\begin{enumerate}[label={\upshape(\roman*)}]
\item\label{Lem:BasicPropertiesi} $\X_1(G)= \Phi(G)$\,;
\item\label{Lem:BasicPropertiesii} $\X_k(G) = \bigcap_{T \in \pc{k-j}(G)} \X_j(T) $ 
for every $j \in [k]$\,;
\item\label{Lem:BasicPropertiesiii} $\X_{i+1}(G) \leq \X_i(G)$\,;
\item\label{Lem:BasicPropertiesiv} $\Phi_i(G) \leq T$, for every $T \in \pc{i}(G)$\,;
\item\label{Lem:BasicPropertiesv} $ \Phi_i(G) \leq \X_i(G)$ for every $i$\,;
\item\label{Lem:BasicPropertiesextra} If $N \unlhd G$
and $N \leq \X_i(G)$ for some index $i$, 
then $\X_i(G/N) = \X_i(G)N/N$.  
\end{enumerate}
\end{lemma}

\begin{proofof}
\ref{Lem:BasicPropertiesi}
Clearly the $\X_i(G)$ are characteristic subgroups of $G$, for all $i$, and 
\[
\X_1(G)= \bigcap_{T \in \pc{1}(G)} T = \Phi(G),
\]
as $\pc{1}(G)$ is exactly the set of maximal subgroups of $G$.

\ref{Lem:BasicPropertiesii} To see this observe that for every $j \in [k]$ we have 
\[
\pc{k}(G)= \bigcup_{T \in \pc{k-j}(G)} \pc{j}(T),
\]
due to the fact that for every subgroup $S$ of index $p^k$ in a $p$-group $G$ 
there exists a subgroup $T\leq G $ of index $p^{k-j}$ that contains $S$. 

\ref{Lem:BasicPropertiesiii} Now this part follows from \ref{Lem:BasicPropertiesii} as 
\[
\X_{i+1}(G) =  \bigcap_{T \in \pc{i}(G)} \X_1(T) 
  \leq \bigcap_{T \in \pc{i}(G)} T 
  =  \X_i(G)\,.
\]

\ref{Lem:BasicPropertiesiv} We induce on $i$. 
For $i = 1$ we clearly have $\Phi(G) \leq M$ 
for every maximal subgroup $M$ of $G$. 
Assuming it holds for $i$, we will prove it for $i+1$. 
Let $T \in \pc{i+1}(G)$ and pick $K \in \pc{i}(G)$ with $T \leq K$. 
So $T$ is a maximal subgroup of $K$ and $\Phi(K) \leq T$. 
Furthermore, the inductive hypothesis implies that 
$\Phi_i(G) \leq K$ and therefore 
\[
\Phi_{i+1}(G) = \Phi(\Phi_i(G)) \leq \Phi(K) \leq T. 
\]

\ref{Lem:BasicPropertiesv} 
The proof of this part follows directly from~\ref{Lem:BasicPropertiesiv} 
and the definition of $\X_i(G)$. 

\ref{Lem:BasicPropertiesextra} 
Since $N \leq \X_i(G)$, it follows from the correspondence theorem
that 
\[
\pc{i}(G/N) = \{H/N : H \in \pc{i}(G)\};
\]
this clearly proves the claim, and the proof of the lemma is complete.

\end{proofof}

\begin{figure}[htbp]
\centering
\begin{tikzpicture}[scale=1.0,sgplattice]
  \node[char] at (4.25,0) (1) {\gn{C1}{C_1}};
  \node[char] at (1.12,0.953) (2) {\gn{C2}{C_2}};
  \node[char] at (7.38,0.953) (3) {\gn{C2}{C_2}};
  \node[char] at (4.25,0.953) (4) {\gn{C2}{C_2 = N}};
  \node[char] at (4.25,2.4) (5) {\gn{C2^2}{C_2^2}};
  \node[norm] at (6.38,2.4) (6) {\gn{C4}{C_4}};
  \node[norm] at (8.38,2.4) (7) {\gn{C4}{C_4}};
  \node at (0.125,2.4) (8) {\gn{C4}{C_4}};
  \node at (2.12,2.4) (9) {\gn{C4}{C_4}};
  \node[norm] at (1.12,3.84) (10) {\gn{C2xC4}{C_2{\times}C_4}};
  \node[norm] at (4.25,3.84) (11) {\gn{C2xC4}{C_2{\times}C_4}};
  \node[char] at (7.38,3.84) (12) {\gn{C2xC4}{C_2{\times}C_4}};
  \node[char] at (4.25,4.8) (13) {\gn{C4:C4}{C_4{\rtimes}C_4}};
  \draw[lin] (1)--(2) (1)--(3) (1)--(4) (2)--(5) (3)--(5) (4)--(5) (3)--(6)
     (3)--(7) (2)--(8) (2)--(9) (5)--(10) (8)--(10) (5)--(11) (9)--(11)
     (5)--(12) (6)--(12) (7)--(12) (12)--(13) (10)--(13) (11)--(13);
  \node[cnj=8] {2};
  \node[cnj=9] {2};
\end{tikzpicture}
\caption[The poset of classes of subgroups of $\mathcal{G}$]{The poset of classes of subgroups of   $\mathcal{G}$.}
\label{figure}
\end{figure}
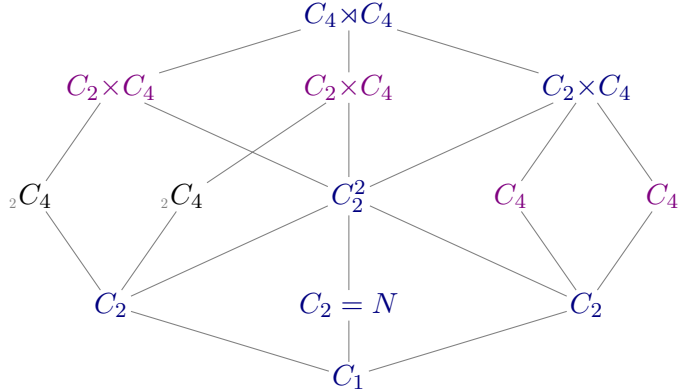

\begin{remark}\label{Rem:FirstRemark}
Firstly, observe that equality really is possible in~\ref{Lem:BasicPropertiesiii},
since in $Q_8$, for example, we have $\X_1(Q_8) = \X_2(Q_8)$.
Secondly, we comment on the behaviour of the $\X$-series 
under quotients. 
Assuming that $G$ is a $p$-group, as usual, 
and that $N$ is an arbitrary normal subgroup of $G$, 
we have $\X_1(G/N) = \X_1(G)N\big/N$
since $\X_1(G) = \Phi(G) = G'\mho(G)$
and each of $G'$, $\mho(G)$ behaves well under quotients.
However, things change for larger-index terms of the $\X$-series.
It is not true, for instance, that $\X_2(G/N) = \X_2(G)N\big/N$
for arbitrary $G$, $N$. 
An example where this fails to be true 
is the group of order $16$ with presentation
\[
\mathcal{G} = C_4 \rtimes C_4 
  = \left\langle a, b\, \big\vert\, a^4 = b^4 = 1, a^b = a^{-1} \right\rangle .
\]
The subgroup $N = \langle a^2b^2\rangle$ is central of order $2$
and affords the quotient $G/N \cong Q_8$.
To see this, observe that $G/N$ is non-abelian
(since, for example, $abN \neq baN$) and 
$a^2N = b^2N$ is the unique element of order $2$ in $G/N$.
Now $\X_2(G)$ is trivial (see Figure~\ref{figure})
and thus $\X_2(G)N\big/N$ is trivial, 
but $\X_2(G/N)$ is not.
\end{remark}

Next, we prove a useful property of the $\X$-series.

\begin{proposition}\label{Prop:SubgMonotone}
Let $G$ be a $p$-group and let $M \leq G$. 
Then $\X_i(M) \leq \X_i(G)$.
\end{proposition}

\begin{proofof}
Fix a subgroup $M \leq G $ and an index $i \geq 0$.
To prove the desired containment, 
it suffices to show that 
$\X_i(M)  \leq K$ for every $K \in \pc{i}(G)$. 
The equality 
% $\frac{|MK|}{|K|}= \frac{|M|}{|M \cap K|}$ 
$|MK|\big/|K| = |M|\big/|M \cap K|$ 
and the fact that $K$ has index $p^i$ in $G$
together imply that 
the index of $K \cap M$ in $M$ is at most $p^i$.
Thus there exists a subgroup $T$ of $K \cap M$ 
whose index in $M$ is exactly $p^i$.
Therefore
\[
\X_i(M)  \leq T \leq K \cap M \leq K.
\]
As this holds for every $K \in \pc{i}(G)$, 
the proposition follows.
\end{proofof}

We will refer to this property of the $\X$-series
as the subgroup-monotonicity of $\X$.
Moreover, the same will be said for a subgroup 
series that satisfies the conclusion of Proposition~\ref{Prop:SubgMonotone} (i.e. that it is subgroup-monotone).

\begin{lemma}\label{Lem:LikeLCS}
Let $G$ be a $p$-group. Then 
\[
\X_i (\X_j(G) ) \leq \X_{i+j} (G).
\]
In addition, $\X_i (G) /\X_{i+1}(G) $ is an elementary abelian $p$-group for every non-negative integer $i$.
\end{lemma}

\begin{proofof}
In view of Lemma~\ref{Lem:BasicProperties} 
we have $\X_{i+j}(G) = \bigcap_{T \in \pc{j}(G)} \X_{i}(T)$.
Let 
\[
K \coloneqq \bigcap_{T \in \pc{j}(G)} T = \X_j(G).
\]
Then $\X_i(K) \leq \X_i(T)$ for every $T \in \pc{j}(G)$ by Proposition~\ref{Prop:SubgMonotone}. 
Hence 
\[
\X_{i+j}(G) = \bigcap_{T \in \pc{j}(G)} \X_{i}(T) \, \geq \, \X_i(K) = \X_i (\X_j(G)).
\]
\noindent
Thus, the first part of the lemma holds, 
from which we deduce that $\X_i (G) /\X_{i+1}(G)$ 
is an elementary abelian group seeing as 
\[
\Phi(\X_i(G)) = \X_1(\X_i(G)) \leq \X_{i+1}(G).
\]
This completes the induction and the proof of the lemma.
\end{proofof}

We now look at the situation where some normal subgroup $N$
of $G$ is complemented in $G$.  Theorem \ref{Thm:converse} is essentially the second part of the following:

\begin{theorem}\label{Thm:DoubleInclusion}
Let $G$ be a  $p$-group and let $N$
be a normal subgroup of $G$.
Suppose that $N$ is complemented, by $H$ say, in $G$.
Then for all $i \geq 0$ we have:
\begin{enumerate}[label={\upshape(\roman*)}]
\item\label{item:1DI} $\X_i(N) \X_i(H) \leq \X_i(G) \leq N \X_i(H)$;
\item\label{item:2DI} $\X_i(H) = \X_i(G) \cap H$.
\end{enumerate}
\end{theorem}

\begin{proofof}
\ref{item:1DI}
The leftmost inclusion is a consequence of the
fact that the product is well-defined (since $\X_i(N)$
is characteristic in $N$ and the latter is normal in $G$)
and the subgroup-monotonicity of the $\X$-series.

As regards the right-hand-side inclusion, 
assume first that $i \leq \log_p(|H|)$.
Observe that for all $T, S \in \pc{i}(H)$
Dedekind's lemma yields
\begin{equation}\label{Eq:1Dedekind}
NS \cap NT = N(NS \cap T)
\end{equation}
and
\begin{equation}\label{Eq:2Dedekind}
NS \cap T = 
NS \cap H \cap T = 
S(N \cap H) \cap T = 
S \cap T.
\end{equation}
By~\eqref{Eq:1Dedekind} and~\eqref{Eq:2Dedekind}
we get $NS \cap NT = N(S \cap T)$ and repeatedly
applying this fact gives us
\[
\bigcap_{T \in \pc{i}(H)} NT = 
N \bigcap_{T \in \pc{i}(H)} T =
N\X_i(H).
\]
But if $T \in \pc{i}(H)$ then $NT \in \pc{i}(G)$,
thus
\[
\X_i(G) =    \bigcap_{S \in \pc{i}(G)} S
        \leq \bigcap_{T \in \pc{i}(H)} NT
        =    N\X_i(H).
\]

In particular, for $i_0 = \log_p(|H|)$ we get
$\X_{i_0}(G) \leq N \X_{i_0}(H) = N$.
Hence for $i \geq i_0$ we see that
$\X_i(G) \leq \X_{i_0}(G) \leq N$.
Therefore the inclusion $\X_i(G) \leq N \X_i(H)$
is valid for all non-negative integers $i$,
completing the proof.

\ref{item:2DI} We have that
\[
\X_i(H) \leq H \cap \X_i(G) \leq H \cap N\X_i(H) =
\X_i(H)(N \cap H) = \X_i(H),
\]
where the first inclusion follows from the 
subgroup-monotonicity of the $\X$-series,
the second inclusion from part~\ref{item:1DI}
and the penultimate equality by Dedekind's lemma.
This double inclusion forces $\X_i(H) = \X_i(G) \cap H$,
as desired.
\end{proofof}

\begin{corollary}\label{Cor:IntersectionX}
Let $G$ be a  $p$-group 
and let $N$ be a normal subgroup of $G$
that  is complemented in $G$ by $H$, say.
Then $N \cap \X_i(G)$ 
has a complement in $\X_i(G)$ 
for all $i \geq 0$
and, in fact, $\X_i(H)$ is one such complement.
\end{corollary}

\begin{proofof}
Fix an index $i \geq 0$.
We show that $\X_i(H)$ is a complement
for $N \cap \X_i(G)$ in $\X_i(G)$
by observing the following:
\[
(N \cap \X_i(G))\X_i(H) = \X_i(G) \cap N\X_i(H) 
                        = \X_i(G).
\]
The first equality is yet another application of Dedekind's
lemma and the second equality follows directly
from Theorem~\ref{Thm:DoubleInclusion}~\ref{item:1DI}. 
The proof is complete.
\end{proofof}

Theorem ~\ref{Thm:DoubleInclusion} will also aid us in the proof
of the next result.

\begin{theorem}\label{Thm:DirectProd}
Given $p$-groups $G$ and $H$, we have 
$\X_i(G \times H) = \X_i(G) \times \X_i(H)$ 
for all non-negative integers $i$. 
\end{theorem}

\begin{proofof}
Fix an index $i$.
Proposition~\ref{Prop:SubgMonotone} clearly implies that 
$\X_i(G) \times 1 = \X_i(G \times 1) \leq \X_i(G \times H)$ 
and similarly for $H$. 
Hence 
\[
\X_i(G) \times \X_i(H) \leq \X_i(G \times H).
\]

\noindent
For the other direction, 
note that by Theorem~\ref{Thm:DoubleInclusion}
we have
\begin{equation}\label{Eq:FirstBound}
\X_i(G \times H) \leq G \X_i(H)
\end{equation}
and 
\begin{equation}\label{Eq:SecondBound}
\X_i(G \times H) \leq \X_i(G) H.
\end{equation}
Combining~\eqref{Eq:FirstBound} and~\eqref{Eq:SecondBound} gives us
\begin{multline*}
\X_i(G \times H)  \leq   G \X_i(H) \cap \X_i(G) H 
= \X_i(G) \left( G \X_i(H) \cap H \right) =  \\
=      \X_i(G) \left[ \X_i(H) (G \cap H) \right]
=      \X_i(G) \times \X_i(H),
\end{multline*}
where all equalities are applications of Dedekind's lemma.
Since both inclusions are valid, 
the proof of the theorem is complete.
\end{proofof}

\begin{proposition}\label{prop:abelian}
Assume that $G$ is an abelian $p$-group. 
Then $\X_i(G) = \Phi_i(G)= \mho_{i}(G)$ and $\X_{j}(\X_{i}(G)) = \X_{i+j}(G)$, 
for all integers $i, j \geq 1$.
\end{proposition}

\begin{proofof}
Assume first that $G$ is a cyclic group. Then 
 $\pc{i}(G) =\mho_{i}(G)$ 
as $\mho_{i}(G)$ is the unique subgroup of $G$ of index $p^i$ in $G$ and thus $\X_i(G)= \mho_{i}(G)= \Phi_i(G)$. 

In the general case, let $G= C_i \times C_2 \times \cdots \times C_t$, 
where $C_j$ are cyclic subgroups of $G$. 
Then in view of Theorem~\ref{Thm:DirectProd} 
we get $\X_i(G) = \prod_{j=1}^t \X_i(C_j)$. 
But for any cyclic subgroup $C$ we have 
$\X_i(C) = \mho_i(C)= \Phi_i(C)$.
Hence 
\[
\X_i(G) = \prod_{j=1}^t \X_i(C_j) 
    = \prod_{j=1}^t \mho_i(C_j)
    = \prod_{j=1}^t \Phi_i(C_j)
    = \Phi_i(G)
    = \mho_i(G),
\]
where the penultimate equality follows from the fact 
that the Frattini subgroup respects direct products 
and an easy induction argument (see~\cite[Satz 6]{gasch}).

\noindent
The last part of the proposition follows directly 
from the corresponding relations that the agemo subgroups 
of an abelian group\footnote{In fact, this
relation among the agemo subgroups is satisfied
more generally by regular $p$-groups; 
see Lemma~1.2.12~(ii) in \cite{LGMK}, for example.} satisfy: 
\[
\mho_{i}(\mho_{j}(G)) = \mho_{i+j}(G).
\]
The proof of the proposition is now complete.
\end{proofof}

\begin{proposition}\label{Prop:CyclicCase}
Let $C$ be a cyclic subgroup of $G$ of order $p^n$. 
Then the following are equivalent:
\begin{enumerate}[label={\upshape(\roman*)}]
\item\label{item:Cyclic1} $C \cap \X_n(G)= 1$;  
\item\label{item:Cyclic2} $C$ has a complement in $G$;  
\item\label{item:Cyclic3} $C \cap \X_i(G)= \X_i(C)$ for every non-negative integer $i$.
\end{enumerate}
\end{proposition}

\begin{proofof}
\ref{item:Cyclic1} $\to$ \ref{item:Cyclic2}
We assume that $C \cap \X_n(G) = 1$  
and wish to establish the existence of a complement 
for $C$ in $G$. 
Indeed, if $C \cap T > 1$ for all subgroups $T$ of index $p^n$ in $G$,
then $C \cap T \geq \Omega_1(C) > 1$,  as $\Omega_1(G)$ is the unique subgroup of order $p$ in $C$. 
Therefore $C \cap \X_n(G) \geq \Omega_1(C) > 1$, 
and this  is clearly a contradiction. 
It follows that there exists 
at least one subgroup $T$ of index $p^n$ in $G$ 
whose intersection with $C$ is trivial. 
That subgroup $T$ serves as the desired complement.

\ref{item:Cyclic2} $\to$ \ref{item:Cyclic3}
Let $T$ be, as in the previous paragraph,
a complement for $C$ in $G$.
As $C$ is a cyclic group, 
we see that $\X_i(C) = \mho_i(C)$ 
is the unique subgroup of order $p^{n-i}$ in $C$. 
Now we pick a subgroup series of $G$ 
\[
T = K_n < K_{n-1} < \ldots < K_1 < G = K_0
\]
with $|K_i/K_{i+1}|= p$ 
and we claim that $K _i = T  \X_i(C)$ 
and $C \cap K_i = \X_i(C)$.
We clearly have $K_i = T  (C \cap K_i)$, 
while 
\[ 
|C \cap K_i| = \frac{|K_i|}{|T|} = \frac{|G|}{p^i}|T| = p^{n-i}.
\]
Since $\X_i(C)$ is the unique subgroup of $C$ of order $p^{n-i}$, 
we deduce that $C \cap K_i = \X_i(C)$, and the claim follows. 
 
Now, $K_i$ has index $p^i$ in $G$ 
and therefore $\X_i(G) \leq K_i$. 
We conclude that 
\[
C \cap \X_i(G) \leq C \cap K_i = \X_i(C) \leq C \cap \X_i(G), 
\]
and thus $\X_i(G) \cap C = \X_i (C)$, as desired.

\ref{item:Cyclic3} $\to$ \ref{item:Cyclic1}
This follows easily from the fact that 
$C \cap \X_i(G)= \X_i(C)$ for all $i \geq 0$ 
upon taking $i = n$ and observing that $\X_n(C) = 1$.
\end{proofof}

Our next result outlines which members of the
lower central series of $G$ are guaranteed to be contained
in some term of the $\X$-series of $G$.
As we will see, our result is essentially best possible.

\begin{proposition}\label{Prop:GammaInX}
If $G$ is a  $p$-group 
and $i$ is a positive integer
then 
\[
\gamma_{p^{i-1} + 1}(G) \leq \X_i(G) 
\]
and the index $p^{i-1} + 1$ is optimal.
\end{proposition}

\begin{proofof}
Let $H$ be an arbitrary subgroup of $G$ of index $p^i$ in $G$. 
The regular action of $G$ on the left cosets of $H$ in $G$
gives rise to a homomorphism $\phi : G \to S_{p^i}$ 
whose kernel $N$ is the core of $H$ in $G$; that is 
\[
N = \bigcap_{g \in G} H^g \leq H.
\]
Let $L_i$ be a Sylow $p$-subgroup of $S_{p^i}$. 
Then $G/N$ is isomorphic to a subgroup of $L_i$. 
Now $L_i$ is an iterated wreath product 
\[
\underbrace{C_p \wr \ldots \wr C_p}_{\text{$i$ terms}}
\]
having order $p^{\frac{p^i-1}{p-1}}$ and class $c = p^{i-1}$
(cf. \cite{kaloujnine}). 
Hence $\gamma_{p^{i-1} + 1}(L_i) = 1$, 
so $\gamma_{p^{i-1} + 1}(G/N) = 1$. 
But the lower central series 
behaves well with respect to quotients 
and therefore 
\[
1 = \gamma_{p^{i-1} + 1}(G/N) = \gamma_{p^{i-1} + 1}(G) N /N.
\]
We conclude that $\gamma_{p^{i-1} + 1}(G) \leq N \leq H $ 
and as $H$ was arbitrary we see that 
\[
\gamma_{p^{i-1} + 1}(G) \leq \bigcap_{H \in \pc{i}(G) } H = \X_i(G).
\]

To see that the index $p^{i-1}+1$ is optimal, 
observe that $\gamma_{p^{i-1}+1}(L_i) = 1$
and the index $j = p^{i-1}+1$ is the least positive
integer such that $\gamma_j(L_i) = 1$
(since $L_i$ has class $p^{i-1}$).
In proof of the second claim of the proposition therefore, 
it will be sufficient
to establish that $\X_i(L_i) = 1$ 
for all positive integers $i$.
We argue by induction on $i$,
noting that the base case $i = 1$ is clearly
valid since $L_1$ has prime order $p$ and thus
its Frattini subgroup is trivial.
Assume validity of the claim for $i-1$.
Note that $L_i = L_{i-1} \wr C_p$ and that the base
of $L_i$ 
\[
\underbrace{B = L_{i-1} \times \ldots \times L_{i-1}}_{\text{$p$ factors}}
\]
is a maximal subgroup of $L_i$.

On the other hand, we have 
by Lemma~\ref{Lem:BasicProperties}~\ref{Lem:BasicPropertiesii}
(with $k = i$ and $j = i-1$)
\[
\X_i(L_i) = \bigcap_{M \in \pc{1}(L_i)} \X_{i-1}(M) \leq \X_{i-1}(B).
\]
But 
\[
\X_{i-1}(B) = \X_{i-1}(L_{i-1}) \times \ldots \times \X_{i-1}(L_{i-1}) = 1,
\]
where the first equality follows from Theorem~\ref{Thm:DirectProd} 
while the second equality
is valid by the inductive hypothesis.
Thus, $\X_i(L_i) = 1$, completing the induction 
and the proof.
\end{proofof}

\begin{corollary}\label{Cor:2groups}
Let $G$ be a  $2$-group. 
Then 
$[G, \Phi(G) ] \leq \X_2(G)$.
\end{corollary}

\begin{proofof}
We begin with the observation that
\[
[G, \Phi(G)] = 
[G, G'\mho(G)] = 
[G,G'] [G, \mho(G)] =
\gamma_3(G) [G, \mho(G)].
\]
We assert that both terms in the right-hand-side
of the equality above are contained in $\X_2(G)$.
That $\gamma_3(G) \leq \X_2(G)$ follows from
Proposition~\ref{Prop:GammaInX} with $p = i = 2$.
Thus, it will suffice to prove that 
$[G, \mho(G)] \leq \X_2(G)$.
According to Corollary 1.28 in \cite{LGMK},
we have
\[
[G,\mho(G)] \leq  \mho([G, G])[G, G, G] 
      =    \mho(G') \gamma_3(G)
      \leq  \mho(G') \X_2(G)
\]
when $G$ is a  $2$-group.
In addition, 
\[
\mho(G')    \leq \mho(\Phi(G)) 
            \leq \Phi(\Phi(G)) 
            =    \Phi_2(G) 
            \leq \X_2(G),
\]
where the first containment is justified 
by the subgroup-monotonicity of the agemo subgroup 
(cf. Lemma 1.2.7 in \cite{LGMK}) 
and the last by Lemma~\ref{Lem:BasicProperties}~\ref{Lem:BasicPropertiesv}.
Hence the claim has been proved.
\end{proofof}

The content of Corollary~\ref{Cor:2groups} is no longer
valid if $p > 2$. 
Indeed, for any odd prime $p$ the Sylow $p$-subgroup
$C_p \wr C_p$ of $S_{p^2}$ satisfies 
$\X_2(C_p \wr C_p) = 1$ since it has a maximal
subgroup which is elementary abelian (in particular, its base), 
but its Frattini subgroup is not central.
To see why this last claim is true, assume that $\Phi(G) \leq Z(G)$
and recall that by~\cite{kaloujnine} the class of 
$C_p \wr C_p$ is $p$. Thus, if $\Phi(G) \leq Z(G)$
then $\Phi(G) = Z(G)$, since $Z(G)$ has order $p$.
On the other hand $G' \leq \Phi(G) = Z(G)$, forcing
$\gamma_3(G)$ to be trivial, which is absurd (even for $p=3$).

\section{Complements and abelian subgroups}\label{Sec:3}

If $A$ is an abelian group then a condition
which is both necessary and sufficient 
for a subgroup $B$ of $A$ to have a complement in $A$
exists and is based on Pr\"{u}fer's notion of a \textbf{pure subgroup};
see for example~\cite[Chap. 5]{fuchs}.
We define this next.

\begin{definition}
Let $A$ be an abelian group and let $B$ be a subgroup of $A$. 
Then $B$ is called a pure subgroup of $A$ provided that 
if $a^n \in B$ for $a \in A$ and $n \in \mathbb{N}$ 
then there exists $b \in B$ such that $a^n = b^n$. 
Equivalently, $B$ is pure in $A$ 
if and only if $B^n = B \cap A^n$ 
for all $n \in \mathbb{N}$.
\end{definition}

For the convenience of the reader we present
here an independent proof of the following result.
Note that no originality is claimed here. 
The result itself should be
well-known to abelian group theorists 
(and, in fact, holds more generally than we have stated),
but perhaps not so much to  group theorists in general.

\begin{theorem}\label{thm:pure}
Let $A$ be a  abelian group and $B \leq A$. Then $B$ has a complement in $A$ if and only if $B$ is a pure subgroup of $A$.
\end{theorem}
\begin{proofof}
Suppose first that $B$ is complemented in $A$ and let $C$ be a complement. Then $A$ is the (internal) direct product $A=B \times C$. Thus, for each $a \in A$ there exist unique $b \in B$ and $c \in C$ such that $a=bc$. Therefore, if $a^n \in B$ for some $n \in \mathbb{N}$ then $b^n c^n \in B$ thus $c^n \in B \cap C = 1$. So $a^n = b^n$ and thus $B$ is pure in $A$.

For the converse, suppose that $B$ is pure in $A$ and that $B \leq C \leq A$, where $C/B$ is cyclic. Let $n = |C : B|$ and write $\langle cB \rangle =C/B$. As $c^n \in B$, there exists $b \in B$ such that $c^n = b^n$. Setting $g\coloneqq cb^{-1}$, we have $g^n=1$. Since $gB= cB$, we see that $n=o(gB)$ in $C/B$. It follows that $C = B \times \langle g \rangle $ and thus $B$ is a direct factor of $C$.

Now for the general case, write $A/B$ as a direct product of
cyclic groups $C_i/B$. Let $X_i$ be a complement for $B$ in
$C_i$ and note that $|X_i|= |C_i:B|$. Write $D = \prod_i X_i$.

We have $BD \geq BX_i = C_i$ for all $i$, and thus
$BD = \prod_i C_i = A$. Also,
\[
|D| \leq \prod_i |X_i| = \prod_i |C_i:B| = |A/B|,
\]
and thus 
$A = B \times D$. Our proof is complete.
\end{proofof}

Now we focus our attention on  $p$-groups. 
In view of Gasch\"{u}tz's theorem 
mentioned in Section~\ref{Sec:1}, this is natural.
In this setting, the previous theorem asserts the
following:

Let $G$ be an abelian  $p$-group and let $A$
be a subgroup of $G$. Then $A$ has a complement in $G$
if and only if 
\[
A^n = A \cap G^n \tag*{($\dagger$)}
\]
for all $n \in \mathbb{N}$.

However, the only relevant exponents $n$ in this case are the powers of $p$,
thus $(\dagger)$ is essentially equivalent to 
\[
\mho_i(A) = A \cap \mho_i(G) \tag*{($\ddag$)}
\]
for all non-negative integers $i$.
On the other hand, the $\mho_i$ subgroup
of an  abelian $p$-group is its $\X_i$ subgroup, according to  Proposition \ref{prop:abelian}.
Thus we have the following.

\begin{corollary}
Let $G$ be a  abelian $p$-group 
and let $A$ be a subgroup of $G$. 
Then $A$ has a complement in $G$
if and only if 
\[
A \cap \X_i(G) = \X_i(A)
\]
for all indices $i \geq 0$.
\end{corollary}

We are now ready to prove that
the condition above, which is both necessary
and sufficient for complementation of a subgroup
of an abelian $p$-group, is sufficient in general.

\begin{proofofmain}
We begin with an easy observation; 
it will suffice to establish the claim 
for all indices $i$ such that $p^i$
is at most the exponent of the subgroup $A$.
Moreover, notice that if $A$ is cyclic,
then the equivalence between \ref{item:Cyclic2} and \ref{item:Cyclic3}
in Proposition~\ref{Prop:CyclicCase}
ensures that $A$ is complemented in $G$.

Now we argue by induction on $|A|$, 
noting that we have
addressed the base case of the induction 
in the previous paragraph. 
We may also assume that $A$ is not cyclic. 
So let $C$ be a cyclic subgroup of maximal order in $A$ 
and note that $|C| = \mathrm{exp}(A) = p^r$, say.
By assumption, $A \cap \X_r(G) = \X_r(A) = 1$. 
Then $C \cap \X_r(G) \leq A \cap \X_r(G) = 1$, 
thus $C \cap \X_r(G) = 1$
and the equivalence between \ref{item:Cyclic1} and \ref{item:Cyclic2}
in Proposition~\ref{Prop:CyclicCase}
suffices to guarantee that $C$ has a complement, say $H$, in $G$.

By Dedekind's lemma, $A = C \times B$, 
where $B \coloneqq A \cap H$.
Observe that since $B$ is complemented by $C$ in $A$ and taking
into account Theorem \ref{Thm:converse}, we have $B \cap \X_i(A) = \X_i(B)$
for all $i \leq t$, where $p^t = \mathrm{exp}(B)$.

Since the $\X$-series is subgroup-monotone, we have for
each $i \leq t$ the following:
\[
B \cap \X_i(H) \leq B \cap \X_i(G) = 
B \cap (A \cap \X_i(G)) = B \cap \X_i(A) = \X_i(B).
\]
Also, again by subgroup-monotonicity we have
$\X_i(B) \leq B$ and $\X_i(B) \leq \X_i(H)$ for all appropriate indices $i$.
Thus, $\X_i(B) \leq B \cap \X_i(H)$ and since the reverse containment
also holds we have $\X_i(B) = B \cap \X_i(H)$.
The induction hypothesis applied to $B$ relative to the subgroup $H$
now ensures that there exists a subgroup $D$ in $H$ which complements $B$. 
Then
\[
A \cap D = (A \cap H) \cap D = B \cap D = 1. 
\]
Moreover,
\[
|A| = |C||B| = |G : H||H : D| = |G : D|.
\]
We conclude that $D$ complements $A$ in $G$ thus our induction is complete. 
\end{proofofmain}

We follow up with a few remarks on the content
and proof of Theorem~\ref{Thm:main}.
\begin{remark}
It should be noted that Theorem~\ref{Thm:main}
works for a class of $p$-groups which is strictly
larger than the class of abelian $p$-groups.
For $A$ to have a complement in $G$ it is sufficient,
as the proof of the theorem indicates,
that $A$ have the following property:
\begin{quote}
for every $H \leq A$ there exists a normal (in $H$) 
cyclic  subgroup $1< C \unlhd  H$ that admits a complement in $H$,
\end{quote}
For example, the property above is enjoyed by the dihedral $2$-groups $D_{2^n}$ which are non-abelian for $n \geq 3$.

On the other hand,  Theorem~\ref{Thm:main} is not true unconditionally; some assumption about $A$ \emph{must} be made. 
Take $\mathcal{H} = C_4 \times \mathcal{G}$, where $\mathcal{G}$ is as in Remark~\ref{Rem:FirstRemark}.
Then there exists a maximal subgroup $M$ of $\mathcal{H}$ which satisfies  the condition of Theorem~\ref{Thm:main}, i.e. $\Phi(M) = \Phi(\mathcal{H})$, yet $M$ has no complement in $\mathcal{H}$ since $\Omega_1(\mathcal{H}) = \Phi(\mathcal{H})$.
\end{remark}

\begin{remark}
On the surface, it looks as though a literal 
interpretation of purity might be enough to
guarantee the existence of a complement in
Theorem~\ref{Thm:main}. That is to say,
under the assumptions of the theorem,
perhaps it is true that an abelian subgroup $A$
of $G$ satisfying $A \cap \mho_i(G) = \mho_i(A)$
for all $i \geq 0$ is necessarily complemented in $G$.
However, this is not so;
purity \emph{really must} be stated in an equivalent form.
To see why, consider an odd prime $p$ and let $G$
be a non-abelian $p$-group of exponent $p$.
Then it is easy to see that 
the equality $A \cap \mho_i(G) = \mho_i(A)$
is valid for every (abelian) subgroup
of $G$, but clearly not every abelian subgroup
can be complemented; simply consider a subgroup of
$\X_1(G) = \Phi(G)$ of order $p$.
\end{remark}

\begin{corollary}\label{Cor:central}
Let $G$ be a $p$-group and let $A$ be a central
subgroup of $G$. Then $A$ has a complement in $G$
if and only if $A \cap \X_i(G) = \X_i(A)$
for all $i \geq 0$. 
\end{corollary}

\begin{proofof}
Since $A$ is central in $G$, it is abelian. 
Thus Theorem~\ref{Thm:main} implies that $A$ is complemented in $G$.
Conversely, suppose that $A$ has a complement in $G$, say $H$. Then $H$ is normalised by both $A$ and $H$, thus $G = A \times H$ and $A$ is, in fact, normally complemented in $G$.
Now Theorem~\ref{Thm:converse} implies that
$\X_i(A) = A \cap \X_i(G)$, as desired.
\end{proofof}

In fact, looking carefully at the proof 
of Theorem~\ref{Thm:main} and
extracting the abstract properties of the $\X$-series
used in the proof allows us to state a generalisation
of Theorem~\ref{Thm:main}.
We need a definition first.

\begin{definition}
Given a  $p$-group $G$ we call the series
\[
\W_i(G) = \X_i(G) \cap \B_i(G), 
\]
a \textbf{good series} for $G$ if the following three conditions are satisfied:
\begin{enumerate}[label={\upshape(\roman*)}]
\item $\{\B_i(G)\}$ is a descending normal series of $G$; 
\item $\B_i(G)$ (and thus $\W_i(G)$ as well) is subgroup-monotone;
\item $\B_i(G)$ (and thus $\W_i(G)$ as well) respects direct products;
\item If $C$ is a cyclic group, then $\Omega_1(C) \leq \B_1(C)$ (and  thus  $\Omega_1(C) \leq \W_1(C)$ as well).
\end{enumerate}
\end{definition}

Clearly, $\X_i(G)$ itself is a good series for $G$. 
Other examples of good series are the following:
\begin{enumerate}[label={\upshape(\arabic*)}]
\item $\W_i(G) = \X_i(G) \cap \Omega_1(G)$;
\item $\W_i(G) = \X_i(G) \cap \mho_{i-1}(G)$;
\item $\W_i(G) = \X_i(G) \cap \mho_{i-1}(G) \cap \Omega_1(G)$.
\end{enumerate}

The generalisation asserted earlier now reads as follows.

\begin{theorem}\label{Thm:mainw}
Let $A$ be an abelian subgroup of the  $p$-group
$G$ and suppose that $\W_i$ is a good series for $G$ so that  
\[
A \cap \W_i(G) = \W_i(A)
\]
for all integers $i \geq 1$. 
Then $A$ has a complement in $G$. 
\end{theorem}

Since the proof is nearly identical to that of Theorem~\ref{Thm:main},
we have elected to omit it.

We remark that the following result,
although resembling Corollary~\ref{Cor:IntersectionX},
is ultimately independent, since no normality
assumption is made.

\begin{proposition}
Let $G$ be a  $p$-group 
and $A$ an abelian subgroup of $G$ such that
\[
A \cap \X_i(G) = \X_i(A)
\]
for all non-negative integers $i$.
Then $A \cap \X_i(G) $ has a complement in $\X_i(G)$ 
for all $i \geq 0$.
\end{proposition}

\begin{proofof}
By Theorem~\ref{Thm:main}
we have already established the case $i = 0$, 
as $\X_0(G) = G$ and $A$ has a complement in $G$. 
According to the same theorem, 
the abelian group $A \cap \X_i(G)$ 
has a complement in $\X_i(G)$ provided we can show that 
\begin{equation}\label{eq:AIntersectSeries}
(A \cap \X_i(G)) \cap \X_t(\X_i(G)) = \X_t(A \cap \X_i(G))    
\end{equation}
for every integer $t \geq 1$. 

\noindent 
One inclusion is straightforward, 
as the $\X$-series is subgroup-monotone and therefore 
\[
\X_t(A \cap \X_i(G)) \leq 
(A \cap \X_i(G)) \cap \X_t(\X_i(G)).
\]
For the other inclusion, 
we first observe that 
\[
\X_t(A \cap \X_i(G)) = \X_t(\X_i(A)) = \X_{t+i}(A),
\]
where the last equality follows from Proposition~\ref{prop:abelian}.
Hence 
\[
(A \cap \X_i(G)) \cap \X_t(\X_i(G)) =    A \cap \X_t(\X_i(G))  
                                    \leq A \cap \X_{t+i}(G)    
                                    =    \X_{t+i}(A)           
                                    =    \X_t( A \cap \X_i(G)). 
\]
So we have shown that equation~\eqref{eq:AIntersectSeries} 
is valid for every $t \geq 1$, and the proposition follows.
\end{proofof}

\begin{corollary}\label{Cor:XDecomposition}
Let $G$ be a  $p$-group and let $A$
be a normal abelian subgroup of $G$ with the property
\[
A \cap \X_i(G) = \X_i(A)
\]
for all non-negative integers $i$.
Then $A$ has a complement in $G$
and for every complement $H$ 
the analogous property
\[
H \cap \X_i(G) = \X_i(H)
\]
holds for all $i \geq 0$. As a consequence
$\X_i(G) = \X_i(A) \rtimes \X_i(H)$
for all $i \geq 0$.
\end{corollary}

\begin{proofof}
Since $A \cap \X_i(G) = \X_i(A)$ for all $i \geq 0$,
Theorem~\ref{Thm:main} implies that $A$
is complemented in $G$.
If $H$ is any such complement for $A$,
then the fact that 
$H \cap \X_i(G) = \X_i(H)$
holds for all $i \geq 0$ is a consequence of
Theorem~\ref{Thm:converse}.
Moreover, Corollary~\ref{Cor:IntersectionX} implies
that $A \cap \X_i(G)$,
which by assumption is  $\X_i(A)$, is complemented
in $\X_i(G)$ by $\X_i(H)$ and this holds, of course,
for all $i \geq 0$.
\end{proofof}

We conclude our work by pointing out the limitations
of Theorem~\ref{Thm:main}.

First of all, by virtue of Proposition~\ref{Prop:CyclicCase}
it is the case that if $A$ is cyclic then it is complemented in $G$
if and only if $A \cap \X_i(G) = \X_i(A)$ for all $i \geq 0$.
Moreover, the condition works as \enquote{if and only if}
in case $A$ is central in $G$
owing to Corollary~\ref{Cor:central}.

On the other hand, observe that an abelian subgroup $A$ 
of a $p$-group $G$ can be complemented in $G$
without satisfying the condition 
$A \cap \X_i(G) = \X_i(A)$ for all $i \geq 0$.
Concrete examples abound;
one possibility is to take $A$ to be an elementary
abelian group containing the Frattini subgroup.
If $A$ is complemented in $G$ (and this happens for
many groups, e.g. $D_8$ with $A$ either of the two
non-cyclic maximal subgroups),
where $G$ is non-abelian,
then it cannot possibly be the case that
$A \cap \X_i(G) = \X_i(A)$ for all $i \geq 0$.
Things go awry as quickly as possible,
i.e. for $i = 1$.

In fact, an elementary abelian subgroup $A$
satisfies $A \cap \X_i(G) = \X_i(A)$ for all $i \geq 0$
if and only if $A \cap \Phi(G) = 1$,
thus Theorem~\ref{Thm:main} is only useful when $A$
avoids the Frattini subgroup of $G$.
However, in that case Theorem~\ref{Thm:main} tells us nothing new
since Satz 7 in Gasch\"{u}tz's paper~\cite{gasch} guarantees
that $A$ has a complement in $G$ regardless of whether
$G$ is a $p$-group or not, as long as $G$ is finite.

What is perhaps not so obvious is the following.

\begin{proposition}\label{Prop:NormalComplement}
Let $G$ be a  $p$-group and let $A$
be a normal subgroup of $G$ such that $A \cap \Phi(G) = 1$.
Then $A$ has a normal complement in $G$ and is thus a direct
factor of $G$.
\end{proposition}

\begin{proofof}
First of all, observe that we do not need to assume
that $A$ is abelian since it follows from
\[
\Phi(A) \leq A \cap \Phi(G) = 1
\]
that $A$ must, in fact, be elementary abelian.

Now we argue by induction on $|A|$. The claim is
evident if $A = 1$, so the induction begins and we
assume that $A > 1$. Since $A \cap \Phi(G) = 1$ and
$A > 1$, there exists a maximal subgroup $M$ of $G$ such
that $A$ is not contained in $M$. Write $D \coloneqq A \cap M$.
In case $D = 1$, we are done because $M$ is normal in $G$.
So suppose that $D > 1$. Observe that since $AM > M$
and $M$ is maximal in $G$, we have $G = AM$. Now both
$M$ and $A$ normalise $D$, so $D$ is normal in $G$.

The induction hypothesis applied to $D$ in $G$ ensures
the existence of a normal complement $K$. In other
words, $G = D \times K$. By Dedekind's lemma, we have
$A = D \times (A \cap K)$. The subgroup $A \cap K$ is normal
in $K$ and proper in $A$. 
Also, $\Phi(G) = \Phi(D) \times \Phi(K) = \Phi(K)$.
This holds because $A$ is elementary abelian hence so is $D$.

Now
\[
(A \cap K) \cap \Phi(K) = 
A \cap \Phi(K) = A \cap \Phi(G) = 1. 
\]
Applying the 
induction hypothesis to $A \cap K$ in $K$, we get
$K = (A \cap K) \times H$ for some subgroup $H$.
But now $H$ is normalised by $K$ and also by $D$, 
because $D$ actually centralises it, and since $G = DK$ we see
that $H$ is normal in $G$. 
It is easy to see that $A \cap H = 1$,
so $H$ is the required normal complement and that
completes our induction.
\end{proofof}

It is clear that the same proof would work just as well
if we had assumed more generally that $G$ is a nilpotent
group instead of a $p$-group. 
We preferred to state the proposition in a more restricted
setting to remain in the geist of the rest of our paper.

\begin{ack}
That Proposition~\ref{Prop:NormalComplement} is true
was communicated to the first author, without proof,
by Martin Isaacs during the winter of 2020.
Both authors are grateful to him for this contribution
and for various suggestions to improve the exposition.
Thanks also go to Benjamin Sambale for spotting typos
in an early draft of the paper and for suggesting
a simplification of the proof of Proposition~\ref{Prop:SubgMonotone}.
\end{ack}
\bibliographystyle{amsalpha}

\begin{thebibliography}{LGM02}

\bibitem[FT63]{FTOdd}
W.~Feit and J.~G. Thompson, \emph{Solvability of groups of odd order}, Pacific
  J. Math. \textbf{13} (1963), 775--1029.

\bibitem[Fuc15]{fuchs}
L.~Fuchs, \emph{Abelian {G}roups}, Springer Monogr. Math., Cham: Springer,
  2015.

\bibitem[Gas52]{Gaschutz:1952}
W.~Gasch{\"u}tz, \emph{Zur {E}rweiterungstheorie der endlichen {G}ruppen}, J.
  Reine Angew. Math. \textbf{190} (1952), 93--107.

\bibitem[Gas53]{gasch}
W.~Gasch\"{u}tz, \emph{\"{U}ber die {$\Phi$}-{U}ntergruppe endlicher
  {G}ruppen}, Math. Z. \textbf{58} (1953), 160--170.

\bibitem[Kal46]{kaloujnine}
L.~Kaloujnine, \emph{Sur les {{\(p\)}}-groupes de sylow du groupe
  sym{\'e}trique du degr{\'e} {{\(p^m\)}}. ({Suite} centrale ascendante et
  descendante)}, C. R. Acad. Sci., Paris \textbf{223} (1946), 703--705
  (French).

\bibitem[Kir17]{kirtland}
J.~Kirtland, \emph{Complementation of {N}ormal {S}ubgroups. {I}n {F}inite
  {G}roups}, Berlin: De Gruyter, 2017.

\bibitem[LGM02]{LGMK}
C.~R. Leedham-Green and S.~McKay, \emph{The {S}tructure of {G}roups of {P}rime
  {P}ower {O}rder}, Lond. Math. Soc. Monogr., New Ser., vol.~27, Oxford: Oxford
  University Press, 2002.

\bibitem[Sam22]{sambale}
B.~Sambale, \emph{On the converse of {G}asch\"{u}tz’ complement theorem},
  \url{https://www.iazd.uni-hannover.de/fileadmin/iazd/sambale/pdfs/GaConverse.pdf},
  2022.
\end{thebibliography}

\end{document}